# DICKSON INVARIANTS, REGULARITY AND COMPUTATION IN GROUP COHOMOLOGY

DAVE BENSON

ABSTRACT. In this paper, we investigate the commutative algebra of the cohomology ring $H^*(G,k)$ of a finite group $G$ over a field $k$. We relate the concept of quasi-regular sequence, introduced by Benson and Carlson, to the local cohomology of the cohomology ring. We give some slightly strengthened versions of quasi-regularity, and relate one of them to Castelnuovo–Mumford regularity. We prove that the existence of a quasi-regular sequence in either the original sense or the strengthened ones is true if and only if the Dickson invariants form a quasi-regular sequence in the same sense. The proof involves the notion of virtual projectivity, introduced by Carlson, Peng and Wheeler.

As a by-product of this investigation, we give a new proof of the Bourguiba–Zarati theorem on depth and Dickson invariants, in the context of finite group cohomology, without using the machinery of unstable modules over the Steenrod algebra.

Finally, we describe an improvement of Carlson's algorithm for computing the cohomology of a finite group using a finite initial segment of a projective resolution of the trivial module. In contrast to Carlson's algorithm, ours does not depend on verifying any conjectures during the course of the calculation, and is always guaranteed to work.

## 1. INTRODUCTION

The purpose of this paper is to investigate the cohomology ring $H^*(G,k)$ of a finite group $G$ over a field $k$ of characteristic $p$ from the point of view of commutative algebra. In particular, we formulate the following conjecture, which is one of the central focuses of the paper.

**Conjecture 1.1.** Let $G$ be a finite group, and let $k$ be a field of characteristic $p$. Then the Castelnuovo–Mumford regularity of the cohomology ring is equal to zero:

$$\mathsf{Reg}\, H^*(G,k) = 0.$$

The author is partly supported by a grant from the NSF. During the course of the preparation of this paper, he was a guest of the Mathematical Sciences Research Institute in Berkeley, California.

The approach is a synthesis of many recent threads in representation theory and cohomology of finite groups, and is inspired in part by ideas of Carlson [13], Greenlees [21], and Okuyama and Sasaki [22]. Some background material for this paper may be found in the survey paper [6]. In particular, we refer the reader to the last section of that paper for examples, and generalizations of the above conjecture to the context of the cohomology of compact Lie groups, orientable virtual Poincaré duality groups, and $p$-adic Lie groups.

Our main theorems are as follows. Recall from §2 of [9] that $H^*(G,k)$ always contains a (possibly dilated) copy of the Dickson invariants as a homogeneous set of parameters. See also §8 of this paper, and Wilkerson's primer [26]. The following theorem is a consequence of Theorems 8.5 and 9.6.

**Theorem 1.2.** *The Dickson invariants form a filter-regular system of parameters in $H^*(G,k)$.*

This allows us to give a new proof (see Theorem 9.9) of the Bourguiba–Zarati theorem [10] relating the depth to the Dickson invariants, in the particular case of finite group cohomology, without using the machinery of unstable algebras over the Steenrod algebra.

**Theorem 1.3** (Bourguiba–Zarati)**.** *Let $G$ be a finite group and let $k$ be a field of characteristic $p$. Suppose that $H^*(G,k)$ has depth $d$. Then the first $d$ Dickson invariants, in ascending order of degree, form a regular sequence in $H^*(G,k)$.*

In the language of quasi-regular sequences, introduced in [8] and further investigated in [12, 13], the following is a consequence of Corollary 4.7 and Theorems 8.5 and 9.6.

**Theorem 1.4.** *Let $G$ be a finite group and $k$ be a field of characteristic $p$. Then the following are equivalent.*

  (i) *There exists a quasi-regular sequence in $H^*(G,k)$ (see Definition 3.2).*

  (ii) *The Dickson invariants form a quasi-regular sequence in $H^*(G,k)$.*

  (iii) *Every homogeneous system of parameters which satisfies the rank-restriction condition (see Definition 8.3) is a quasi-regular sequence.*

  (iv) *Every filter-regular sequence in $H^*(G,k)$ is a quasi-regular sequence.*

  (v) *For $0 \le i \le r$, we have $a^i_{\mathfrak{m}} H^*(G,k) < 0$.*

The proof gives something stronger, which we formulate using strengthened versions of the condition of quasi-regularity. In the language of Definition 3.2, the above theorem is still true with quasi-regular replaced by strongly quasi-regular or by very strongly quasi-regular, but with stronger bounds in condition (v). Namely, for strongly quasi-regular, the inequality in (v) is replaced by $a_\mathfrak{m}^i H^*(G,k) \leq -i$ for $0 \leq i \leq r$, and for very strongly quasi-regular, it is replaced by $a_\mathfrak{m}^i H^*(G,k) < -i$ for $0 \leq i \leq r-1$ and $a_\mathfrak{m}^r H^*(G,k) \leq -r$.

We also prove that the equivalent conditions of even the very strongly quasi-regular version always hold if the difference between the Krull dimension and the depth of $H^*(G,k)$ is at most two. This improves a theorem of Okuyama and Sasaki [22] for quasi-regular sequences. The proof of the following theorem can be found in §9.

**Theorem 1.5.** *If the Krull dimension and the depth of $H^*(G,k)$ differ by at most two, then there is a very strongly quasi-regular sequence in $H^*(G,k)$. In particular,* $\mathsf{Reg}\, H^*(G,k) = 0$.

In §10, we use the machinery we have developed to give an improvement of Carlson's method [12, 13] for computing the cohomology of a finite group using a finite initial segment of a projective resolution of the trivial module. Carlson's algorithm depends on some conjectures in the cohomology of finite groups. The conjectures are verified for the group in question during the course of the computation, so there is no uncertainty about the answer. But the method is not guaranteed to succeed, if one of the conjectures should fail for the group in question. Our improved version avoids the use of these conjectures, and gives improved bounds on how much must be computed before the computation can be deemed to be complete.

Theorem 10.1 is a slightly strengthened version of the following statement. Suppose that we've computed up to degree $N$, and that it's at least far enough so that we've at least got a homogeneous sequence of parameters $\zeta_1, \ldots, \zeta_r$ in $H^*(G,k)$ (this can be checked by restriction to elementary abelian $p$-subgroups). Consider the ring $\tau_N H^*(G,k)$ defined by the generators and relations we've found so far. If $\zeta_1, \ldots, \zeta_r$ are filter-regular in $\tau_N H^*(G,k)$ and

$$N > \max\{0, \mathsf{Reg}\, \tau_N H^*(G,k)\} + \sum_{j=1}^r (|\zeta_j| - 1) \qquad (1.6)$$

then the obvious map $\tau_N H^*(G,k) \to H^*(G,k)$ is an isomorphism. The given condition is guaranteed to hold for large enough $N$, and in the presence of Conjecture 1.1, the bound (1.6) becomes $N > \sum_{j=1}^r (|\zeta_j| - 1)$.

The inequality $>$ in (1.6) can be replaced by $\geq$ if the $p$-rank of the center of a Sylow $p$-subgroup of $G$ is at least two, or if we have some other *a priori* way of knowing that the depth of the cohomology ring is at least two.

**Acknowledgement.** I would like to thank Rodney Sharp for his comments on an earlier draft of this paper.

2. Notation and conventions

We shall need to use the theory of varieties for modules, so we begin by introducing the notation. For background, see for example [1, 20]. Let $G$ be a finite group, and $k$ be a field of characteristic $p$. We write $V_G$ for the maximal ideal spectrum of $H^*(G,k)$. This is a homogeneous affine variety whose dimension is equal to the $p$-rank $r_p(G)$, which we shall denote by $r$.

If $M$ is a finitely generated $kG$-module, we write $V_G(M)$ for the subvariety of $V_G$ defined by the kernel of the action of $H^*(G,k)$ on $\operatorname{Ext}^*_{kG}(M,M)$. Then

$$V_G(M_1 \oplus M_2) = V_G(M_1) \cup V_G(M_2)$$
$$V_G(M_1 \otimes M_2) = V_G(M_1) \cap V_G(M_2).$$

A $kG$-module $M$ is projective if and only if $V_G(M) = \{0\}$. If $H$ is a subgroup of $G$, the restriction map $\operatorname{res}_{G,H} \colon H^*(G,k) \to H^*(H,k)$ gives rise to a map $\operatorname{res}^*_{G,H} \colon V_H \to V_G$.

If $\zeta$ is a nonzero element of $H^n(G,k)$, then we write $L_\zeta$ for the kernel of a representative cocycle $\hat\zeta \colon \Omega^n k \to k$. This gives a short exact sequence

$$0 \to L_\zeta \to \Omega^n k \to k \to 0. \tag{2.1}$$

The variety $V_G(L_\zeta)$ is equal to the hypersurface $V_G\langle \zeta \rangle$ in $V_G$ defined by regarding $\zeta$ as an element of the coordinate ring. If $\zeta_1, \ldots, \zeta_r$ is a homogeneous set of parameters in $H^*(G,k)$ then

$$\begin{aligned} V_G(L_{\zeta_1} \otimes \cdots \otimes L_{\zeta_r}) &= V_G(L_{\zeta_1}) \cap \cdots \cap V_G(L_{\zeta_r}) \\ &= V_G\langle\zeta_1\rangle \cap \cdots \cap V_G\langle\zeta_r\rangle = \{0\} \end{aligned}$$

and so $L_{\zeta_1} \otimes \cdots \otimes L_{\zeta_r}$ is projective.

Similarly, the variety of $L_{\zeta_1} \otimes \cdots \otimes L_{\zeta_{r-1}}$ is one dimensional. So this module is periodic.

We work throughout this paper with graded rings satisfying the following conditions.

**Condition 2.2.** Let $k$ be a field of characteristic $p$ (possibly $p = 0$). The $\mathbb{Z}$-graded ring $H$ is finitely generated by positive degree elements as an algebra over $H^0 = k$, and is graded commutative in the sense that $ab = (-1)^{|a||b|}ba$. The maximal ideal $\mathfrak{m}$ is generated by the positive degree elements, so that

$$H = H^0 \oplus H^{>0} = k \oplus \mathfrak{m}.$$

If $M$ is a graded module over $H$, we use the notation $M[i]$ for the shifted module with $M[i]^n = M^{i+n}$. The element of $M[i]$ corresponding to $m \in M$ is written $\sigma^i m$. The action of $x \in H^j$ is given by $x\,\sigma^i m = (-1)^{ij}\sigma^i(xm)$.

## 3. Filter-regular sequences

We shall need the following condition on a sequence of integers.

**Condition 3.1.** The sequence of integers $d_0, \ldots, d_r$ satisfies

(a) $d_i \geq d_{i-1} - 1$ for $1 \leq i \leq r$,

(b) $d_{i-1} \geq d_i$ for $1 \leq i \leq r$, and

(c) $d_0 \geq -1$.

**Definition 3.2.** Let $H$ be a graded ring satisfying Condition 2.2, let $M$ be a finitely generated graded $H$-module of Krull dimension $r$, and let $d_0, \ldots, d_r$ be integers satisfying Condition 3.1(a) and (b). Let $\zeta_1, \ldots, \zeta_r$ be a homogeneous system of parameters in $H$ for $M$ (i.e., $M/(\zeta_1, \ldots, \zeta_r)$ has finite length) with $|\zeta_i| = n_i$. We say that $\zeta_1, \ldots, \zeta_r$ is a *filter-regular* system of parameters of type $(d_0, \ldots, d_r)$ for $M$ if for each $i = 0, \ldots, r-1$, the map

$$(M/(\zeta_1, \ldots, \zeta_i))^j \to (M/(\zeta_1, \ldots, \zeta_i))^{j+n_{i+1}} \tag{3.3}$$

induced by multiplication by $\zeta_{i+1}$ is injective for $j > n_1 + \cdots + n_i + d_i$, and $(M/(\zeta_1, \ldots, \zeta_r))^j = 0$ for $j > n_1 + \cdots + n_r + d_r$. If there exist values of $d_0, \ldots, d_r$ for which this condition holds, then we say that $\zeta_1, \ldots, \zeta_r$ is *filter-regular* for $M$. The existence of a filter-regular sequence is guaranteed by the standard method of prime avoidance. The smaller the values of the $d_i$, the stronger the condition becomes. If $M = H$ and $d_i = -1$ for $0 \leq i \leq r$ then the sequence is *quasi-regular*, in the terminology of [8].

We introduce two strengthened versions of the above quasi-regularity condition. If $d_i = -i$ for $0 \leq i \leq r$, then we say that the sequence is *strongly quasi-regular*. We shall see that the existence of a strongly quasi-regular sequence is equivalent to $\mathsf{Reg}\, H = 0$, so that Conjecture 1.1 states that the cohomology of a finite group has a strongly quasi-regular sequence.

If $d_i = -i - 1$ for $0 \leq i \leq r - 1$ and $d_r = -r$, then we say that the sequence is *very strongly quasi-regular*. This is the strongest of these conditions for which the method of §7 works. There is quite a lot of evidence that the cohomology ring of a finite group always has a very strongly quasi-regular sequence. We shall prove this in the case where the Krull dimension and depth differ by at most two.

**Lemma 3.4.** *Suppose that $M$ has depth $d$. Then the first $d$ terms of any filter-regular sequence form a regular sequence.*

*Proof.* (See also exercise 18.3.8 of Brodmann and Sharp [11]). Let $\zeta_1, \ldots, \zeta_r$ be a filter-regular sequence. We prove by induction on $i$ that for $0 \leq i \leq d - 1$, multiplication by $\zeta_{i+1}$ is injective on $M/(\zeta_1, \ldots, \zeta_i)$. It follows from the inductive hypothesis that $M/(\zeta_1,\ldots,\zeta_i)$ has depth $d-i > 0$. The existence of a regular element in this quotient means that if multiplication by $\zeta_{i+1}$ has a kernel then the kernel has elements of arbitrarily large degree. So $\zeta_{i+1}$ is regular, and the proof by induction is complete. □

Filter-regularity can be interpreted in a number of cohomological ways, and this is the subject of the next few sections. We begin by interpreting it in terms of the Koszul complex. We write $\mathcal{K}^*(M;\zeta_1,\ldots,\zeta_i)$ for the doubly graded Koszul complex of $M$ with respect to parameters $\zeta_1,\ldots,\zeta_i$. In order to fit with the indexing of the spectral sequences of §6, we regard the degree coming from $M$ as vertical, and we make the Koszul degree horizontal and negative. So $\mathcal{K}^0$ is $M$ with no shift, and $\mathcal{K}^{-i}$ is $M$ shifted upwards in degree by $\sum_{j=1}^{i}|\zeta_j|$.

**Lemma 3.5.** *Let $d_0,\ldots,d_r$ be integers satisfying Condition 3.1(a) and (b). If $\zeta_1,\ldots,\zeta_r$ is a filter-regular sequence of type $(d_0,\ldots,d_r)$ for $M$ with $|\zeta_i|=n_i$ then*

$$H^{-s,t}(\mathcal{K}^*(M;\zeta_1,\ldots,\zeta_r)) = 0 \quad \text{for} \quad t > n_1 + \cdots + n_s + d_{r-s}.$$

*Proof.* Write $_iH^{-s,t}$ for $H^{-s,t}(\mathcal{K}^*(M;\zeta_1,\ldots,\zeta_i))$. Then the short exact sequence of complexes

$$0 \to \mathcal{K}^*(M;\zeta_1,\ldots,\zeta_{i-1}) \to \mathcal{K}^*(M;\zeta_1,\ldots,\zeta_i) \to \mathcal{K}^{*+1}(M;\zeta_1,\ldots,\zeta_{i-1})[-n_i] \to 0$$

gives rise to a long exact sequence

$$\cdots \to {_{i-1}}H^{-s,t} \to {_i}H^{-s,t} \to {_{i-1}}H^{-s+1,t-n_i} \xrightarrow{\zeta_i} {_{i-1}}H^{-s+1,t} \to \cdots$$

This allows us to show by induction on $i$, using the inequality $d_i \geq d_{i-1} - 1$, that $_iH^{-s,t} = 0$ provided $s > 0$ and $t > n_1 + \cdots + n_s + d_{i-s}$. If $s = 0$ then $_rH^{0,t} = (M/(\zeta_1,\ldots,\zeta_r))^t$ is zero for $t > n_1 + \cdots + n_r + d_r$. □

Sometimes it is more convenient to work in terms of the stable Koszul complex, or in other words, in terms of local cohomology. This is the subject of the next section.

## 4. Local cohomology

Recall that for graded modules over a graded ring, local cohomology has two gradings. We write $H_\mathfrak{m}^{i,n}M$ for the degree $n$ part of the $i$th local cohomology module of $M$. So the first grading gives local cohomological degree, and the second keeps track of internal degree.

The definition of Castelnuovo–Mumford regularity is usually only considered in the context of strictly commutative graded rings generated in degree one. The following definition from Eisenbud and Goto [18] seems to work well in our more general context.

**Definition 4.1.** Let $H$ be a graded ring satisfying Condition 2.2. If $M$ is a finitely generated graded $H$-module, we set
$$a_\mathfrak{m}^i(M) = \max\{n \in \mathbb{Z} \mid H_\mathfrak{m}^{i,n} M \neq 0\}$$
(or $a_\mathfrak{m}^i(M) = -\infty$ if $H_\mathfrak{m}^i M = 0$). The *Castelnuovo–Mumford regularity* of $M$ is defined to be
$$\mathsf{Reg}\, M = \max_{j \geq 0}\{a_\mathfrak{m}^j(M) + j\}.$$

Remarks concerning the history and geometric significance of this definition can be found in §20.5 of Eisenbud [17].

**Theorem 4.2.** *Let $G$ be a finite group of $p$-rank $r$ and let $k$ be a field of characteristic $p$. Then $\mathsf{Reg}\, H^*(G, k) \geq 0$.*

*Proof.* The "last survivor" described in theorem 1.3 of [8] and reinterpreted in terms of local cohomology in theorem 4.1 of [5] (see also §5 below) says that $H_\mathfrak{m}^{r,-r} H^*(G, k) \neq 0$. So $a_\mathfrak{m}^r(H^*(G, k)) \geq -r$. □

The proof of the following lemma is similar to arguments from Trung [24, 25].

**Lemma 4.3.** *Let $H$ be a graded ring satisfying Condition 2.2, and let $M$ be a finitely generated graded $H$-module. If $\zeta \in H^n$ is such that the kernel of $\zeta$ on $M$ is $\mathfrak{m}$-torsion, then*
$$a_\mathfrak{m}^{i+1}(M) + n \leq a_\mathfrak{m}^i(M/\zeta M) \leq \max(a_\mathfrak{m}^i(M), a_\mathfrak{m}^{i+1}(M) + n). \quad (4.4)$$

*Proof.* The short exact sequence
$$0 \to \mathrm{Ker}(\zeta) \to M \to M/\mathrm{Ker}(\zeta) \to 0$$
shows that for $i > 0$ we have $H_\mathfrak{m}^i M = H_\mathfrak{m}^i(M/\mathrm{Ker}(\zeta))$. So for $i \geq 0$ the exact sequence
$$0 \to M/\mathrm{Ker}(\zeta)[-n] \xrightarrow{\zeta} M \to M/\zeta M \to 0$$
gives a long exact sequence
$$\cdots \to H_\mathfrak{m}^{i,j} M \to H_\mathfrak{m}^{i,j}(M/\zeta M) \to H_\mathfrak{m}^{i+1,j-n} M \xrightarrow{\zeta} H_\mathfrak{m}^{i+1,j} M \to \cdots$$
The element of highest degree in $H_\mathfrak{m}^{i+1} M$ is in the kernel of $\zeta$, so this gives the desired inequalities. □

**Theorem 4.5.** *Let $H$ be a graded ring satisfying Condition 2.2, let $M$ be a finitely generated graded $H$-module of Krull dimension $r$, and let $d_0, \ldots, d_r$ be integers satisfying Condition 3.1(a) and (b). Then the following are equivalent.*
  (i) *There exists a filter-regular sequence $\zeta_1, \ldots, \zeta_r$ of type $(d_0, \ldots, d_r)$ for $M$.*
  (ii) *Every filter-regular sequence $\zeta_1, \ldots, \zeta_r$ for $M$ has type $(d_0, \ldots, d_r)$.*
  (iii) *For $0 \leq i \leq r$ we have $a_\mathfrak{m}^i(M) \leq d_i$.*

*If these equivalent conditions hold, then*
$$\text{Reg}\, M \le \max_{0 \le i \le r}\{d_i + i\}.$$

*Proof.* (i) ⇒ (iii): Suppose that $\zeta_1, \ldots, \zeta_r$ is a filter-regular sequence of type $(d_0, \ldots, d_r)$ for $M$ with $|\zeta_i| = n_i$. We use downward induction on $s$, beginning with the case $s = r$, to show that for all $i \ge 0$ and $0 \le s \le r$ we have
$$a_{\mathfrak{m}}^i(M/(\zeta_1, \ldots, \zeta_s)) \le n_1 + \cdots + n_s + d_{s+i}. \tag{4.6}$$
The case $i = 0$ comes from the definition of filter-regular sequence of type $(d_0, \ldots, d_s)$. The case $i > 0$ comes from the inductive hypothesis and the left half of the inequality (4.4). The case $s = 0$ gives (iii).

(iii) ⇒ (ii): If $a_{\mathfrak{m}}^i(M) \le d_i$ for all $i \ge 0$ and $\zeta_1, \ldots, \zeta_r$ is filter-regular, then we use (upward) induction on $s$ to prove that (4.6) holds for all $i \ge 0$ and all $0 \le s \le r$. The case $s = 0$ is the given condition. For $s > 0$, using the right half of inequality (4.4), we have
$$a_{\mathfrak{m}}^i(M/(\zeta_1, \ldots, \zeta_s)) \le \max(a_{\mathfrak{m}}^i(M/(\zeta_1, \ldots, \zeta_{s-1})), a_{\mathfrak{m}}^{i+1}(M/(\zeta_1, \ldots, \zeta_{s-1})) + n_s)$$
$$\le \max(n_1 + \cdots + n_{s-1} + d_{s+i-1}, n_1 + \cdots + n_s + d_{s+i}).$$
Condition 3.1(a) ensures that
$$n_1 + \cdots + n_{s-1} + d_{s+i-1} \le n_1 + \cdots + n_s + d_{s+i}$$
since $n_s \ge 1$. So by induction we obtain inequality (4.6) for all $i \ge 0$ and $0 \le s \le r$. The case $i = 0$ gives (ii).

The implication (ii) ⇒ (i) is obvious. Finally, if (iii) holds then the last statement of the theorem is immediate from the definition of regularity. □

We spell out several particular cases of this theorem.

**Corollary 4.7.** *Suppose that $H$ is a graded ring satisfying Condition (2.2). Then*

(a) *The following are equivalent.*
  (i) *There is a quasi-regular sequence for $H$,*
  (ii) *Every filter-regular sequence for $H$ is quasi-regular,*
  (iii) *For all $i \ge 0$ we have $a_{\mathfrak{m}}^i(H) < 0$.*
*Under these conditions we have $\text{Reg}\, H < r$.*

(b) *The following are equivalent.*
  (i) *There is a strongly quasi-regular sequence for $H$,*
  (ii) *Every filter-regular sequence for $H$ is strongly quasi-regular,*
  (iii) $\text{Reg}\, H \le 0$ *(i.e., for all $i \ge 0$ we have $a_{\mathfrak{m}}^i(H) \le -i$).*

(c) *The following are equivalent.*

(i) *There exists a very strongly quasi-regular sequence for $H$,*
(ii) *Every filter-regular sequence for $H$ is very strongly quasi-regular,*
(iii) *For $0 \le i \le r-1$ we have $a_\mathfrak{m}^i(H) < -i$, and $a_\mathfrak{m}^r(H) \le -r$.*

*Under these conditions, we have $\mathsf{Reg}\, H \le 0$.* □

**Corollary 4.8.** *There exists a strongly quasi-regular sequence for $H^*(G,k)$ if and only if $\mathsf{Reg}\, H^*(G,k) = 0$.*

*Proof.* This follows from Theorem 4.2 and Corollary 4.7. □

## 5. Free resolutions

In this section we relate the invariants $a_\mathfrak{m}^i(M)$ to some invariants defined in terms of free resolutions over polynomial algebras. The theorems of this section generalize some of the theorems of Eisenbud and Goto [18], with proofs similar to the outline given in exercise 15.3.7 of Brodmann and Sharp [11].

Suppose that $H$ satisfies Condition 2.2, and that $M$ is a finitely generated graded $H$-module of Krull dimension $r$. Using Noether normalization, we choose a polynomial ring $R = k[\zeta_1, \ldots, \zeta_r]$ with $|\zeta_i| = n_i$ and a homomorphism of graded rings $R \to H$ so that $M$ is a finitely generated graded $R$-module, and $\zeta_1, \ldots, \zeta_r$ is filter-regular for $M$. Let

$$0 \to F_r \to \cdots \to F_0 \to M \to 0$$

be the minimal resolution of $M$ over $R$. Then we define $\beta_j^R(M)$ to be the largest degree of a generator of $F_j$ as an $R$-module (or $\beta_j^R(M) = -\infty$ if $F_j = 0$). This is the same as the largest value of $\ell$ such that $\mathrm{Tor}_{j,\ell}^R(M,k) \ne 0$.

**Theorem 5.1.** $\beta_0^R(M) \le \max_i\{a_\mathfrak{m}^i(M) + i\} + \sum_{j=1}^r (n_j - 1)$.

*Proof.* We prove by downward induction on $s$ that

$$\beta_0^R(M/(\zeta_1,\ldots,\zeta_s)) \le \max_i\{a_\mathfrak{m}^i(M/(\zeta_1,\ldots,\zeta_s)) + i\} + \sum_{j=s+1}^r (n_j - 1). \quad (5.2)$$

This is clearly true for $s = r$, so suppose $s < r$. By the left hand inequality of (4.4), we have

$$a_\mathfrak{m}^{i+1}(M/(\zeta_1,\ldots,\zeta_s)) + i + 1 + \sum_{j=s+1}^r (n_j - 1)$$
$$\le a_\mathfrak{m}^i(M/(\zeta_1,\ldots,\zeta_{s+1})) + i + \sum_{j=s+2}^r (n_j - 1).$$

Now given a minimal set of generators for $M/(\zeta_1,\ldots,\zeta_{s+1})$, we can lift to elements of $M/(\zeta_1,\ldots,\zeta_s)$ and adjoin a set of $\mathfrak{m}$-torsion elements to obtain a set of generators for $M/(\zeta_1,\ldots,\zeta_s)$. So we have

$$\beta_0^R(M/(\zeta_1,\ldots,\zeta_s)) \le \max\{a_\mathfrak{m}^0(M/(\zeta_1,\ldots,\zeta_s)), \beta_0^R(M/(\zeta_1,\ldots,\zeta_{s+1}))\}$$

$$\leq \max\{a_{\mathfrak{m}}^0(M/(\zeta_1,\ldots,\zeta_s)) + \sum_{j=s+1}^r (n_j-1), \beta_0^R(M/(\zeta_1,\ldots,\zeta_{s+1}))\}$$

So if (5.2) is true for $s+1$ then it is true for $s$. $\square$

**Corollary 5.3.** *Let $d_0,\ldots,d_r$ be integers satisfying Condition 3.1(a). Then*
$$\beta_0^R(M) - d_r \leq \max_i\{a_{\mathfrak{m}}^i(M) - d_i\} + \sum_{j=1}^r n_j.$$

*Proof.* This follows from the theorem together with the inequality
$$-d_r \leq -d_i - i + r,$$
which in turn follows from Condition 3.1(a). $\square$

**Lemma 5.4.** *Let $F$ be a finitely generated free graded $R$-module (i.e., a finite direct sum of shifts of copies of $R$). Then $a_{\mathfrak{m}}^r(F) = \beta_0^R(F) - \sum_{j=1}^r n_j$.*

*Proof.* A direct computation using the Koszul complex shows that $a_{\mathfrak{m}}^r(R) = -\sum_{j=1}^r n_j$. The module $F$ consists of a direct sum of shifts of copies of $R$, with top generator in degree $\beta_0^R(F)$. $\square$

The following is the main theorem of this section.

**Theorem 5.5.** *Let $d_0,\ldots,d_r$ be integers satisfying Condition 3.1(a) and (b), and suppose that $M$, $H$ and $R$ are as described at the beginning of this section. Then*
$$\max_{i\geq 0}\{a_{\mathfrak{m}}^i(M) - d_i\} = \max_{i\geq 0}\{\beta_{r-i}^R(M) - d_i\} - \sum_{j=1}^r n_j. \tag{5.6}$$

*Proof.* We prove this by induction on the projective dimension of $M$ as an $R$-module. In the case of projective dimension zero, the theorem follows directly from Lemma 5.4. So suppose that $M$ has projective dimension at least one. Let $0 \to N \to F \to M \to 0$ be the first stage of a minimal resolution of $M$ over $R$. Set $d_i' = d_{i-1}$ if $i > 0$ and $d_0' = d_0$. Then the sequence of integers $d_0',\ldots,d_r'$ satisfies Condition 3.1(a) and (b). We have $H_{\mathfrak{m}}^i F = 0$ for $i \neq r$, and $a_{\mathfrak{m}}^r(F) = \beta_0^R(M) - \sum_{j=0}^r n_j$ by Lemma 5.4. So the long exact sequence in local cohomology
$$\cdots \to H_{\mathfrak{m}}^{i-1} N \to H_{\mathfrak{m}}^{i-1} F \to H_{\mathfrak{m}}^{i-1} M \to H_{\mathfrak{m}}^i N \to H_{\mathfrak{m}}^i F \ldots$$
shows that for $0 < i < r$ we have $a_{\mathfrak{m}}^i(N) = a_{\mathfrak{m}}^{i-1}(M)$, so that
$$a_{\mathfrak{m}}^i(N) - d_i' = a_{\mathfrak{m}}^{i-1}(M) - d_{i-1},$$
while $a_{\mathfrak{m}}^0(N) = -\infty$. Furthermore,
$$a_{\mathfrak{m}}^{r-1}(M) \leq a_{\mathfrak{m}}^r(N) \leq \max\{a_{\mathfrak{m}}^{r-1}(M), \beta_0^R(M) - \sum_{j=1}^r n_j\},$$
so that using Condition 3.1(b), we have
$$a_{\mathfrak{m}}^{r-1}(M) - d_{r-1} \leq a_{\mathfrak{m}}^r(N) - d_r'$$

$$\leq \max\{a_{\mathfrak{m}}^{r-1}(M) - d_{r-1}, \beta_0^R(M) - d_{r-1} - \sum_{j=1}^r n_j\}$$
$$\leq \max\{a_{\mathfrak{m}}^{r-1}(M) - d_{r-1}, \beta_0^R(M) - d_r - \sum_{j=1}^r n_j\},$$

and

$$a_{\mathfrak{m}}^r(M) \leq a_{\mathfrak{m}}^r(F) = \beta_0^R(M) - \sum_{j=1}^r n_j.$$

On the other hand, for $0 < i \leq r$ we have $\beta_{r-i}^R(N) = \beta_{r-i+1}^R(M)$, so that

$$\beta_{r-i}^R(N) - d_i' = \beta_{r-i+1}^R(M) - d_{i-1},$$

and $\beta_r^R(N) = -\infty$.

Combining these statements, we obtain

$$\max_{i\geq 0}\{a_{\mathfrak{m}}^i(M) - d_i\} \leq \max\{\max_{i\geq 0}\{a_{\mathfrak{m}}^i(N) - d_i'\}, a_{\mathfrak{m}}^r(M) - d_r\}$$
$$\leq \max\{\max_{i\geq 0}\{a_{\mathfrak{m}}^i(N) - d_i'\}, \beta_0^R(M) - d_r - \sum_{j=1}^r n_j\}.$$

Combining this with Corollary 5.3, we obtain

$$\max_{i\geq 0}\{a_{\mathfrak{m}}^i(M) - d_i\} = \max\{\max_{i\geq 0}\{a_{\mathfrak{m}}^i(N) - d_i'\}, \beta_0^R(M) - d_r - \sum_{j=1}^r n_j\}.$$

By the inductive hypothesis, the theorem is true for $N$. So we have

$$\max_{i\geq 0}\{a_{\mathfrak{m}}^i(M) - d_i\} = \max\{\max_{i\geq 0}\{\beta_{r-i}^R(N) - d_i'\}, \beta_0^R(M) - d_r\} - \sum_{j=1}^r n_j$$
$$= \max_{i\geq 0}\{\beta_{r-1}^R(M) - d_i\} - \sum_{j=1}^r n_j.$$

This completes the inductive proof of the theorem. □

Note that the summation term in part (ii) of the following corollary is usually omitted in the literature because of the assumption that the generators are in degree one.

**Corollary 5.7.** *We have*
(i) $\max_{i\geq 0}\{a_{\mathfrak{m}}^i(M)\} = \max_{i\geq 0}\{\beta_i^R(M)\} - \sum_{j=1}^r n_j,$
(ii) $\operatorname{Reg} M = \max_{i\geq 0}\{a_{\mathfrak{m}}^i(M) + i\} = \max_{i\geq 0}\{\beta_i^R(M) - i\} - \sum_{j=1}^r (n_j - 1).$ □

**Corollary 5.8.** (i) *There exists a quasi-regular sequence in $H$ if and only if for all $i \geq 0$ we have $\beta_i^R(H) < \sum_{j=1}^r n_j$.*

(ii) *There exists a strongly quasi-regular sequence in $H$ if and only if for all $i \geq 0$ we have $\beta_i^R(H) \leq i + \sum_{j=1}^r (n_j - 1)$.*

(iii) *There exists a very strongly quasi-regular sequence in $H$ if and only if $\beta_0^R(H) \leq \sum_{j=1}^r (n_j - 1)$, and for $1 \leq i \leq r$, $\beta_i^R(H) < i + \sum_{j=1}^r (n_j - 1)$.* □

## 6. The last survivor

In this section, we construct spectral sequences closely related to those of Benson and Carlson [8], Carlson [12] and Greenlees [21]. We describe the "last survivor" from §7 of [8] in terms of these spectral sequences.

Let $\zeta_1, \ldots, \zeta_r$ be a homogeneous system of parameters in $H^*(G, k)$ with $|\zeta_i| = n_i$. then the last survivor is a nonzero element of $H^*(G, k)$ of degree $\sum_{i=1}^{r}(n_i - 1)$, not in the ideal generated by the parameters, whose existence is guaranteed from the spectral sequence. This discussion will be needed in the proof of Theorem 7.1, and again in §10. The proof given here for the existence of the last survivor is significantly simpler than the one to be found in [8].

Let $L_{\zeta_i}$ be the modules defined by (2.1). Set

$$M_0 = k, \qquad M_i = L_{\zeta_1} \otimes \cdots \otimes L_{\zeta_i} \quad (i \geq 1).$$

Note that $M_r$ is projective, and $M_{r-1}$ is periodic, by the argument given in §2.

Tensoring the short exact sequence (2.1) for $L_{\zeta_i}$ with $M_{i-1}$ gives a short exact sequence

$$0 \to M_i \to \Omega^{n_i} M_{i-1} \oplus \text{(projective)} \to M_{i-1} \to 0. \tag{6.1}$$

To build the spectral sequence, for each parameter $\zeta_i$ we truncate the sequence (2.1) to form a cochain complex

$$0 \to \Omega^{n_i} k \xrightarrow{\hat{\zeta}_i} k \to 0 \tag{6.2}$$

where $k$ is in degree zero and $\Omega^{n_i} k$ is in degree $-1$. Let $C_i^*$ be the tensor product of these complexes for $\zeta_1, \ldots, \zeta_i$. Let $\hat{P}_*$ be a complete resolution of $k$ as a $kG$-module, and set $_iE_0^{-s,t} = \operatorname{Hom}_{kG}(\hat{P}_t, C_i^{-s})$. Then

$$_iE_1^{-s,t} = \hat{H}^t(G, C_i^{-s}) \cong \mathcal{K}^*(\hat{H}^*(G, k); \zeta_1, \ldots, \zeta_i),$$

the Koszul complex for $\hat{H}^*(G, k)$ with respect to the parameters $\zeta_1, \ldots, \zeta_i$. The $d_1$ differential is the same as the Koszul differential. This $_iE_1$ page is indexed so that it lives in columns $-i$ to $0$ in the plane, with each column consisting of a number of copies of Tate cohomology, suitably shifted. The copy in degree zero is unshifted. Using the other spectral sequence of this double complex, we see that the above spectral sequence converges to $\hat{H}^{t-s+i}(G, M_i)$. So we have

$$_iE_1^{**} = \mathcal{K}^*(\hat{H}^*(G, k); \zeta_1, \ldots, \zeta_i) \Rightarrow \hat{H}^{*+i}(G, M_i).$$

In particular, if $i = r$, we have

$$_rE_1^{**} = \mathcal{K}^*(\hat{H}^*(G, k); \zeta_1, \ldots, \zeta_r) \Rightarrow 0. \tag{6.3}$$

This spectral sequence satisfies a duality which combines the Poincaré duality of [8] with Tate duality, as explained in §6 of [7]. We recall that if we dualize the complex (6.2), we obtain a complex of the form

$$0 \to k \xrightarrow{\pm \hat{\zeta}_i} \Omega^{-n_i} k \to 0$$

where the sign is $(-1)^{n_i(n_i+1)/2}$ (see §5 of [8]). In particular, $L^*_{\zeta_i} \cong \Omega^{-n_i-1} L_{\zeta_i}$, $M^*_i \cong \Omega^{-n_1-\cdots-n_i-i} M_i$, and $\hat{H}^j(G, M^*_i) \cong \hat{H}^{j+n_1+\cdots+n_i+i}(G, M_i)$.

Combining this with Tate duality, we obtain a duality between ${}_i E_1^{-s,t}$ and ${}_i E_1^{s-i,n_1+\cdots+n_i-t-1}$. Since this duality is induced by a homomorphism of double complexes at the ${}_i E_0$ level, it commutes with the differentials in the spectral sequence and passes down to a duality on each page. The duality on ${}_i E_\infty^{**}$ is a graded version of the duality between $\hat{H}^{t-s+i}(G, M_i)$ and

$$\hat{H}^{(s-i)+(n_1+\cdots+n_i-t-1)+i}(G, M_i) \cong \hat{H}^{s-t-i-1}(G, M^*_i).$$

If $H^*(G, k)$ is Cohen–Macaulay, then $\zeta_1, \ldots, \zeta_r$ is a regular sequence. In this case, provided $r \geq 2$, there are no nonzero products in Tate cohomology from negative to positive degrees, by lemma 2.1 and theorem 3.1 of [7]. It follows that ${}_r E_2^{**}$ is concentrated in the $-r$ and $0$ columns. So $d_j = 0$ for $2 \leq j < r$, and $d_r$ is an isomorphism from the $-r$ column to the $0$ column.

If $E$ is an elementary abelian $p$-subgroup of $G$ of maximal rank (i.e., rank $r$), then $H^*(E, k)$ is Cohen–Macaulay. Restriction from $G$ to $E$ induces a map of the above spectral sequences. It follows that the identity element in ${}_r E_1^{0,0}$ is not hit by any differential until the ${}_r E_r$ page. So there is a uniquely defined element of ${}_r E_r^{-r,r-1}$ which is sent by $d_r$ to the identity element. The image of this element under the map ${}_r E_r^{-r,r-1} \subseteq {}_r E_1^{-r,r-1} \xrightarrow{\cong} \hat{H}^{-n_1-\cdots-n_r+r-1}(G, k)$ is written $\gamma(\zeta_1, \ldots, \zeta_r)$ (see §7 of [8]). This element, or an element of positive Tate cohomology with nonzero pairing with this element, is called the "last survivor." Note that the element of negative degree is well defined, whereas the element of positive degree is really a coset. It is also observed in [8] that since $\gamma(\zeta_1, \ldots, \zeta_r)$ has nonzero restriction to every elementary abelian $p$-subgroup of maximal rank, and restriction is Tate dual to transfer, the last survivor in positive degree can be taken to be a transfer from any given elementary abelian $p$-subgroup of maximal rank.

Now there are maps of cochain complexes

$$\begin{array}{ccccccccc}
0 & \longrightarrow & 0 & \longrightarrow & k & \longrightarrow & 0 \\
& & \downarrow & & \downarrow & & \\
0 & \longrightarrow & \Omega^{n_i} k & \xrightarrow{\hat{\zeta}_i} & k & \longrightarrow & 0.
\end{array}$$

So if $1 \leq i < r$, then there is a map of complexes $C^*_{i-1} \to C^*_i$ and hence a map of spectral sequences ${}_{i-1} E^{**}_* \to {}_i E^{**}_*$. Composing, we get maps of spectral sequences from each ${}_i E^{**}_*$ to ${}_r E^{**}_*$, which do the obvious thing at the $E_1$ page. Since the identity element in $_rE_1^{0,0}$ is the image under $d_r$ of the element coming from the last survivor, it follows that for $i < r$, the identity element in $_iE_1^{0,0}$ survives to $_iE_\infty^{0,0}$. Applying the duality, $_iE_1^{-i,n_1+\cdots+n_i-1}$ is a one dimensional space isomorphic to $\hat{H}^{-1}(G,k)$, and this one dimensional space survives to $_iE_\infty^{-i,n_1+\cdots+n_i-1}$ provided $i < r$.

## 7. Module theoretic interpretation

In this section, we interpret filter-regularity of a sequence of parameters in $H^*(G,k)$ in module theoretic terms.

**Theorem 7.1.** *Let $d_0,\ldots,d_r$ be integers satisfying Condition 3.1. Then a sequence of parameters $\zeta_1,\ldots,\zeta_r$ for $H^*(G,k)$ is filter-regular of type $(d_0,\ldots,d_r)$ if and only if for each $i = 0,\ldots,r-1$, multiplication by $\zeta_{i+1}$ is injective on $H^j(G,M_i)$ for $j > n_1 + \ldots n_i + i + d_i$.*

*Under these circumstances, if $1 \le i \le r-1$ and $j > n_1 + \cdots + n_i + d_i$, then we have*

$$H^{j+i}(G, M_i) \cong (H^*(G,k)/(\zeta_1,\ldots,\zeta_i))^j.$$

*Proof.* For $j \ge n_1 + \cdots + n_i$, this theorem is an easy inductive argument, using the long exact sequence in cohomology coming from the short exact sequence (6.1). If some of the $d_i$ are negative, the improvement in the bound comes from a careful investigation of Tate cohomology in degree $-1$.

We begin the induction as follows. The long exact sequence in Tate cohomology coming from the short exact sequence (2.1) for $\zeta_1$ is

$$\cdots \to \hat{H}^{j+n_1}(G,M_1) \to \hat{H}^j(G,k) \xrightarrow{\zeta_1} \hat{H}^{j+n_1}(G,k) \to \hat{H}^{j+n_1+1}(G,M_1) \to \cdots$$

By a theorem of Duflot [16], the depth of $H^*(G,k)$ is at least one. So by Lemma 3.4, every filter-regular sequence begins with a regular element. It follows that in the above sequence, the map given by multiplication by $\zeta_1$ is injective for $j \ge 0$. So we have $H^{j+1}(G,M_1) \cong (H^*(G,k)/(\zeta_1))^j$ for $j > n_1 - 1$. Provided $r > 1$, multiplication by $\zeta_1$ is the zero map on $\hat{H}^{-1}(G,k)$, and so $H^{n_1}(G,M_1) \cong H^{n_1-1}(G,k)$ as required.

Next, the long exact sequence in Tate cohomology coming from the short exact sequence (6.1) for $\zeta_2$ is

$$\cdots \to \hat{H}^{j+n_2}(G,M_2) \to \hat{H}^j(G,M_1) \xrightarrow{\zeta_2} \hat{H}^{j+n_2}(G,M_1) \to \hat{H}^{j+n_2+1}(G,M_2) \to \cdots$$

Provided $j \geq n_1$, we have

$$\begin{array}{ccc}
(\hat{H}^*(G,k)/(\zeta_1))^{j-1} & \xrightarrow{\zeta_2} & (\hat{H}^*(G,k)/(\zeta_1))^{j+n_2-1} \\
\cong \downarrow & & \cong \downarrow \\
\hat{H}^j(G, M_1) & \xrightarrow{\zeta_2} & \hat{H}^{j+n_2}(G, M_1)
\end{array}$$

and so $\zeta_2$ is injective on $\hat{H}^j(G, M_1)$ for $j > n_1 + d_1$ if and only if it is injective on $(H^*(G,k)/(\zeta_1))^j$ for $j > n_1 + d_1 - 1$, except possibly in the case where $d_1 = -2$ (its smallest possible value) and $j = n_1 - 1$. We also deduce, for the same values of $j$, that

$$H^{j+n_2+1}(G, M_2) \cong (H^*(G, M_1)/(\zeta_2))^{j+n_2} \cong (H^*(G,k)/(\zeta_1, \zeta_2))^{j+n_2-1}.$$

In the case where $d_1 = -2$ and $j = n_1 - 1$, we must work harder. In this case, we have the following commutative diagram.

$$\begin{array}{ccccccccc}
\hat{H}^{-2}(G,k) & \xrightarrow{\zeta_1} & H^{n_1-2}(G,k) & \longrightarrow & H^{n_1-1}(G, M_1) & \longrightarrow & \hat{H}^{-1}(G,k) & \xrightarrow{0} & \\
\zeta_2 \downarrow & & \zeta_2 \downarrow & & \downarrow d_2 & & \downarrow 0 & & \\
H^{n_2-2}(G,k) & \xrightarrow{\zeta_1} & H^{n_1+n_2-2}(G,k) & \twoheadrightarrow & H^{n_1+n_2-1}(G, M_1) & \xrightarrow{0} & H^{n_2-1}(G,k) & & \\
& & & & \downarrow & & & & \\
& & & & H^{n_1+n_2}(G, M_2) & & & &
\end{array}$$

Now $\hat{H}^{-1}(G,k)$ is one dimensional. The issue is whether an element $x$ of $H^{n_1-1}(G, M_1)$ with nonzero image $\alpha$ in $\hat{H}^{-1}(G,k)$ can make $H^{n_1+n_2}(G, M_2)$ too small. Lift $\zeta_2 x \in H^{n_1+n_2-1}(G, M_1)$ to an element $y$ of $H^{n_1+n_2-2}(G,k)$. This is well defined modulo $(\zeta_1, \zeta_2)$, and is in fact the Massey product $\langle \zeta_1, \alpha, \zeta_2 \rangle$. The crucial observation is that this element is $d_2(\alpha) \in H^*(G,k)/(\zeta_1, \zeta_2)$, where $d_2$ is the differential in the spectral sequence $_2E_2^{**}$ constructed in the last section. Provided $r \geq 3$, we showed that this differential is zero. It follows that $y$ can be adjusted by an element in the image of $\zeta_1$ so that $y = \zeta_2 z$ with $z \in H^{n_1-2}(G,k)$. The image $u$ of $z$ in $H^{n_1-1}(G, M_1)$ has the property that $\zeta_2 u = \zeta_2 x$ in $H^{n_1+n_2-1}(G, M_1)$. This shows that

$$H^{n_1+n_2}(G, M_2) \cong (H^*(G,k)/(\zeta_1, \zeta_2))^{n_1+n_2-2}$$

as required.

The inductive case is similar, but it is clear that it is more efficient to express the last step directly in terms of the spectral sequence. The long exact sequence

$$\cdots \to \hat{H}^{j+n_i}(G, M_i) \to \hat{H}^j(G, M_{i-1}) \xrightarrow{\zeta_i} \hat{H}^{j+n_i}(G, M_{i-1}) \to \hat{H}^{j+n_i+1}(G, M_i) \to \cdots$$

deals with everything, except in the case where $d_i = -i-1$ (its smallest possible value), and then there is a question as to what happens in the lowest degree. In terms of the spectral sequence, the information we have is that $_iE_2^{s,t} = 0$ provided $s < 0$ and $s+t \geq n_1 + \cdots + n_i - i - 1$, except for $_iE_2^{-i,n_1+\cdots+n_i-1}$, which is isomorphic to the one dimensional space $\hat{H}^{-1}(G,k)$. Provided $i < r$, this element survives to $_iE_\infty$ by the discussion at the end of the previous section. The inductive step then follows from the fact that this spectral sequence converges to $\hat{H}^{*+i}(G, M_i)$.

At the end of the induction, notice that the equivalent condition does not mention $d_r$. This is because $M_r$ is a projective module (see §2), and so $H^j(G, M_r) = 0$ for all $j > 0$. So $\zeta_r \colon H^j(G, M_{r-1}) \to H^j(G, M_{r-1})$ is an isomorphism for $j > 0$, and it follows that $(H^*(G,k)/(\zeta_1, \ldots, \zeta_r))^j = 0$ for $j > n_1 + \cdots + n_r + d_{r-1}$. Since $d_r \geq d_{r-1}$, the last condition for filter-regularity holds, and we are done. □

**Corollary 7.2.** *Suppose that $\zeta_1, \ldots, \zeta_r$ is a filter-regular system of parameters of type $(d_0, \ldots, d_r)$ for $H^*(G,k)$. Then $\zeta_1, \ldots, \zeta_r$ is also filter-regular of type $(d_0, \ldots, d_{r-2}, d_{r-2} - 1, d_{r-2} - 1)$.*

*Proof.* As mentioned in §2, the module $M_{r-1}$ is periodic, and $M_r$ is projective. So multiplication by $\zeta_r$ is an isomorphism from $H^j(G, M_{r-1})$ to $H^{j+n_r}(G, M_{r-1})$ for $j > 0$, and $H^j(G, M_r) = 0$ for $j > 0$. It follows that the condition on $d_{r-1}$ in Theorem 7.1 is automatically satisfied. There is no condition on $d_r$ in the theorem except $d_r \geq d_{r-1}$. □

We shall improve on this corollary in Theorem 9.10.

## 8. Dickson invariants

In this section we review the construction of the Dickson invariants in the cohomology ring of a finite group, using the theory of Chern classes. We begin by reminding the reader that the cohomology ring of $BU(n)$ is a polynomial ring on the Chern classes $c_1, \ldots, c_n$,

$$H^*(BU(n); \mathbb{Z}) = \mathbb{Z}[c_1, \ldots, c_n],$$

where $|c_i| = 2i$. If $R$ is a commutative ring of coefficients then the universal coefficient theorem gives $H^*(BU(n); R) = R[c_1, \ldots, c_n]$. The *total Chern class* is defined to be the inhomogeneous element $c = 1 + c_1 + c_2 + \cdots$.

The maximal torus in $U(n)$ is $U(1)^n$, whose cohomology is $R[x_1, \ldots, x_n]$ with $|x_i| = 2$. The restriction map $H^*(BU(n); R) \to H^*(BU(1)^n; R)$ takes $c_i$ to the $i$th elementary symmetric function of $x_1, \ldots, x_n$.

Let $G$ be a finite group. A complex unitary representation $\rho \colon G \to U(n)$ gives rise to a restriction map in cohomology $\rho^* \colon H^*(BU(n); R) \to H^*(BG; R)$.

Recall that the Chern classes of $\rho$ are defined to be $c_i(\rho) = \rho^*(c_i) \in H^*(BG; R) = H^*(G, R)$. For example, if $\rho$ is a one dimensional representation, there is only one Chern class $c_1(\rho)$. It is the image of $\rho$ under the map

$$\mathrm{Hom}(G, U(1)) \cong H^1(G, U(1)) \cong H^2(G, \mathbb{Z}) \to H^2(G, R).$$

The total Chern class of a direct sum of representations is the product of the total Chern classes, $c(\rho_1 \oplus \rho_2) = c(\rho_1)c(\rho_2)$. We are interested in what happens when $\rho = \rho_G$ is the regular representation.

If $E$ is an elementary abelian group of order $p^r$, we write $H^\bullet(G, k)$ for the subring of $H^*(E, k)$ generated by the Bocksteins of the degree one elements. So $H^\bullet(E, k)$ is a polynomial ring on $r$ generators in degree two. The regular representation of $E$ decomposes as a direct sum of all the one dimensional representations, once each. So $\rho_E$ factors as $E \to U(1)^{p^r} \to U(p^r)$. As $\rho$ runs over the one dimensional representations, $c_1(\rho)$ runs over all the elements of

$$H^{2\cdot}(E, \mathbb{F}_p) = H^2(E, \mathbb{F}_p) \cap H^\bullet(E, \mathbb{F}_p) \subseteq H^2(E, k).$$

So $c(\rho_E) = \prod_{w \in H^{2\cdot}(E, \mathbb{F}_p)} (1 + w)$, which means that $c_i(\rho_E)$ is the $i$th elementary symmetric function of the elements of $H^{2\cdot}(E, \mathbb{F}_p)$.

We recall from §8.1 of [2] that the invariants of $\mathrm{Aut}(E) = GL(r, \mathbb{F}_p)$ on $H^\bullet(E, k)$ are given by

$$H^\bullet(E, k)^{GL(r, \mathbb{F}_p)} = k[c_{r,r-1}, \ldots, c_{r,0}]$$

where the Dickson invariants $c_{r,r-j}$ are defined by the equation

$$\prod_{w \in H^{2\cdot}(E, \mathbb{F}_p)} (X - w) = X^{p^r} - c_{r,r-1}X^{p^{r-1}} + c_{r,r-2}X^{p^{r-2}} - \cdots + (-1)^r c_{r,0}X.$$

In particular, we see from this that the $i$th elementary symmetric function of the elements of $H^{2\cdot}(E, \mathbb{F}_p)$ vanishes except when $i$ is of the form $p^r - p^{r-j}$, in which case it is the Dickson invariant $c_{r,r-j}$. Interpreting this in terms of Chern classes, we see that for an elementary abelian $p$-group of order $p^r$, the Chern class $c_i(\rho_E)$ vanishes except when $i$ is of the form $p^r - p^{r-j}$, in which case it is equal to the Dickson invariant $c_{r,r-j} \in H^{2(p^r-p^{r-j})}(E, k)$, so that

$$c(\rho_E) = 1 + c_{r,r-1} + \cdots + c_{r,0}.$$

Now if $S$ is any $p$-group of order $p^n$ and $S'$ is a subgroup of $S$ of index $p^a$, then the restriction of $\rho_S$ to $S'$ is a direct sum of $p^a$ copies of $\rho_{S'}$. So $\mathrm{res}_{S,S'} c(\rho_S) = c(\rho_{S'})^{p^a}$. In particular, if $E$ is an elementary abelian subgroup of order $p^r$ then

$$\mathrm{res}_{S,E} c(\rho_S) = c(\rho_E)^{p^{n-r}} = 1 + c_{r,r-1}^{p^{n-r}} + c_{r,r-2}^{p^{n-r}} + \cdots + c_{r,0}^{p^{n-r}}.$$

So we define $D_j(S) = c_{p^n-p^{n-j}}(\rho_S)$, and we have
$$\operatorname{res}_{S,S'} D_j(S) = D_j(S')^{|S:S'|}, \tag{8.1}$$
and in particular
$$\operatorname{res}_{S,E} D_j(S) = \begin{cases} c_{r,r-j}^{p^{n-r}} & j \leq r \\ 0 & j > r. \end{cases}$$

If $G$ is any finite group, let $S$ be a Sylow $p$-subgroup of $G$, of order $p^n$. Then in the language of the Cartan–Eilenberg stable element method [15], $D_j(S)$ is a stable element in $H^*(S,k)$ for the fusion in $G$. So it is an element of $H^*(G,k) \subseteq H^*(S,k)$. So we write $D_j(G) \in H^{2(p^n-p^{n-j})}(G,k)$ for these elements.

We have proved the following theorem on the existence of Dickson invariants in the cohomology ring of a finite group.

**Theorem 8.2.** *Let $G$ be a finite group and $k$ a field of characteristic $p$. If $|G|_p = p^n$, then there are elements $D_j(G) \in H^{2(p^n-p^{n-j})}(G,k)$ for $j \geq 0$, with the following properties.*

(i) *If $H$ is a subgroup of $G$ then $\operatorname{res}_{G,H} D_j(G) = D_j(H)^{|G:H|_p}$.*

(ii) *If $E$ is an elementary abelian group of order $p^r$ then $D_j(E)$ is the Dickson invariant $c_{r,r-j}$ if $j \leq r$, and $D_j(E) = 0$ if $j > r$.* □

We interpret the $p$-power powers in the exponents in the above theorem as saying that the Dickson invariants are *dilated* by these powers of $p$. Of course, there is often a copy of the Dickson invariants in $H^*(G,k)$ with much smaller dilation. For example, in [4] it is shown that the mod two cohomology of the sporadic Conway group $Co_3$ contains a copy of the rank four Dickson invariants with no dilation, namely a polynomial ring on generators of degrees 8, 12, 14 and 15.

**Definition 8.3.** We say that a homogeneous system of parameters $\zeta_1, \ldots, \zeta_r$ in $H^*(G,k)$ satisfies the *rank-restriction condition* if, for $1 \leq i \leq r$, $\zeta_i$ restricts to zero on every elementary abelian $p$-subgroup of $G$ of rank less than $i$.

**Lemma 8.4.** *Suppose that $\zeta_1 \ldots, \zeta_r$ is a homogeneous set of parameters for $H^*(G,k)$ that satisfies the rank-restriction condition. If $E$ is an elementary abelian $p$-subgroup of $G$ of rank $i$ then the restrictions of $\zeta_1, \ldots, \zeta_i$ to $H^*(E,k)$ form a homogeneous system of parameters, and $\zeta_{i+1}, \ldots, \zeta_r$ restrict to zero.*

*Proof.* It is clear from the definition that $\zeta_{i+1}, \ldots, \zeta_r$ restrict to zero on $H^*(E,k)$. By Evens' finite generation theorem [19], $H^*(E,k)$ is finitely generated as a module over the image of restriction from $G$. Since $H^*(G,k)$ is finitely generated as a module over $k[\zeta_1, \ldots, \zeta_r]$, it follows that $H^*(E,k)$ is finitely generated as a module over the subring generated by the restrictions of $\zeta_1, \ldots, \zeta_i$. Since $H^*(E,k)$ has Krull dimension $i$, it follows that these restrictions form a homogeneous system of parameters for $H^*(E, k)$. $\square$

**Theorem 8.5.** *The Dickson invariants satisfy the rank-restriction condition.*

*Proof.* It follows from Theorem 8.2 that if $E$ is an elementary abelian $p$-subgroup of $G$ of rank less than $i$ then $\operatorname{res}_{G,E}(D_i(G)) = D_i(E)^{|G:E|_p} = 0$. $\square$

We shall make use of Theorem 8.5 in §9.

## 9. Virtual projectivity

The purpose of this section is to prove Theorem 9.6. The proof involves the notion of virtual projectivity introduced by Carlson, Peng and Wheeler [14]. We begin by recalling the basic definitions. Let $W$ be a finitely generated $kG$-module. A sequence
$$0 \to M_1 \to M_2 \to M_3 \to 0 \tag{9.1}$$
is said to be *$W$-split* if the sequence
$$0 \to M_1 \otimes W \to M_2 \otimes W \to M_3 \otimes W \to 0$$
obtained by tensoring with $W$ splits.

A module $M$ is said to be *relatively $W$-projective* if every homomorphism to the module $M_3$ in a $W$-split sequence (9.1) lifts to $M_2$.

For any $kG$-modules $M$ and $N$, there is a transfer map
$$\operatorname{Tr}_W \colon \operatorname{Hom}_{kG}(M \otimes W, N \otimes W) \to \operatorname{Hom}_{kG}(M, N).$$
For example, if $W$ is the permutation module on the cosets of a subgroup $H$ then the image of $\operatorname{Tr}_W$ agrees with the image of $\operatorname{Tr}_{H,G}$. More generally, if $W$ is the direct sum of such permutation modules as $H$ runs over the elements of a collection $\mathcal{H}$ of subgroups of $G$, then the image of $\operatorname{Tr}_W$ is equal to the sum of the images of $\operatorname{Tr}_{H,G}$ for $H \in \mathcal{H}$.

We write $V_G(\operatorname{Tr}_W)$ for the subvariety of $V_G$ defined by the image of
$$\operatorname{Tr}_W \colon \operatorname{Hom}_{kG}(W, W) \to \operatorname{Hom}_{kG}(k, k).$$

**Theorem 9.2.** *Let $W$ be a finitely generated $kG$-module and let $M$ be any $kG$-module. Then the following are equivalent.*

(i) *$M$ is relatively $W$-projective.*

(ii) *$M$ is isomorphic to a direct summand of a tensor product of some module with $W$.*

(iii) *The identity map in $\operatorname{Hom}_{kG}(M, M)$ can be written as $\operatorname{Tr}_W(\gamma)$ for some $\gamma \in \operatorname{Hom}_{kG}(M \otimes W, M \otimes W)$.*

*Proof.* See propositions 2.4 and 3.2 of [14]. $\square$

It follows from this theorem that if $M$ is relatively $W$-projective then every element of $\operatorname{Ext}^*_{kG}(M,M)$ is in the image of $\operatorname{Tr}_W$. A module $M$ is said to be *virtually $W$-projective* if all elements of large enough degree in $\operatorname{Ext}^*_{kG}(M,M)$ are in the image of $\operatorname{Tr}_W$. An example is given in [14] of a module that is virtually $W$-projective but not relatively $W$-projective.

**Theorem 9.3.** *Let $W$ and $M$ be finitely generated $kG$-modules such that*
$$V_G(M) \cap V_G(\operatorname{Tr}_W) = \{0\}.$$
*Then $M$ is virtually $W$-projective.*

*Proof.* See proposition 5.3 of [14]. □

Now let $\mathcal{H}_i$ be the collection of centralizers $C_G(E)$ of elementary abelian $p$-subgroups $E$ of rank $i$. Let $W_i$ be the direct sum of the permutation modules corresponding to the subgroups in $\mathcal{H}_i$, so that the image of $\operatorname{Tr}_{W_i}$ is the sum of the images of $\operatorname{Tr}_{H,G}$ as $H$ runs over the elements of $\mathcal{H}_i$.

**Lemma 9.4.** *The subvariety $V_G(\operatorname{Tr}_{W_i})$ is contained in the union of the images of $\operatorname{res}^*_{E,G}\colon V_E \to V_G$ as $E$ runs over the set of elementary abelian $p$-subgroups of $G$ of rank strictly less than $i$.*

*Proof.* This follows from the transfer theorem of [3]. □

**Theorem 9.5.** *Suppose that $\zeta_1, \ldots, \zeta_r$ is a homogeneous system of parameters for $H^*(G,k)$ that satisfies the rank-restriction condition. Then*
$$M_{i-1} = L_{\zeta_1} \otimes \cdots \otimes L_{\zeta_{i-1}}$$
*is virtually $W_i$-projective.*

*Proof.* By Lemmas 8.4 and 9.4, we have
$$V_G(M_{i-1}) \cap V_G(\operatorname{Tr}_{W_i}) = V_G\langle\zeta_1\rangle \cap \cdots \cap V_G\langle\zeta_{i-1}\rangle \cap V_G(\operatorname{Tr}_{W_i})$$
$$\subseteq V_G\langle\zeta_1\rangle \cap \cdots \cap V_G\langle\zeta_{i-1}\rangle \cap \bigcup_{r_p(E)<i} \operatorname{Im}(\operatorname{res}^*_{G,E}) = \{0\}.$$

It now follows from Theorem 9.3 that $M_{i-1}$ is virtually $W_i$-projective. □

**Theorem 9.6.** *If $\zeta_1, \ldots, \zeta_r$ satisfy the rank-restriction condition, then they form a filter-regular sequence.*

*Proof.* By Duflot's theorem [16], if $H = C_G(E) \in \mathcal{H}_i$ then the restrictions of $\zeta_1, \ldots, \zeta_{i-1}$ form a regular sequence in $H^*(H,k)$. So by Theorem 7.1, multiplication by $\zeta_{i-1}$ is injective on $H^j(H, M_{i-2})$ for $j$ large enough.

Now consider the diagram

$$\begin{array}{ccccccc}
\cdots \longrightarrow & H^{j+n_{i-1}}(G, M_{i-1}) & \longrightarrow & H^{j}(G, M_{i-2}) & \xrightarrow{\zeta_{i-1}} & H^{j+n_{i-1}}(G, M_{i-2}) \\
& \uparrow \sum \mathrm{Tr}_{H,G} & & \uparrow \sum \mathrm{Tr}_{H,G} & & \uparrow \sum \mathrm{Tr}_{H,G} \\
\cdots \longrightarrow & \bigoplus_{H \in \mathcal{H}_i} H^{j+n_{i-1}}(H, M_{i-1}) & \longrightarrow & \bigoplus_{H \in \mathcal{H}_i} H^{j}(H, M_{i-2}) & \xrightarrow{\zeta_{i-1}} & \bigoplus_{H \in \mathcal{H}_i} H^{j+n_{i-1}}(H, M_{i-2})
\end{array} \tag{9.7}$$

We have just seen that the map marked $\zeta_{i-1}$ on the bottom row is injective for all large enough $j$. By Theorem 9.5, the left hand vertical map is surjective for all large enough $j$. A diagram chase shows that the map marked $\zeta_{i-1}$ on the top row is injective for all large enough $j$. $\square$

**Corollary 9.8.** *Let $G$ be a finite group of $p$-rank $r$ and let $k$ be a field of characteristic $p$. Then the Dickson invariants $D_1, \ldots, D_r$ form a filter-regular sequence in $H^*(G, k)$.*

*Proof.* This follows from Theorems 8.5 and 9.6. $\square$

The following is the Bourguiba–Zarati theorem in the particular case of finite group cohomology. Their proof uses the classification of injective unstable modules over the Steenrod algebra.

**Theorem 9.9.** *Let $G$ be a finite group and let $k$ be a field of characteristic $p$. Suppose that $H^*(G, k)$ has depth $d$. Then the Dickson invariants $D_1, \ldots, D_d$ form a regular sequence in $H^*(G, k)$.*

*Proof.* This follows from Lemma 3.4 and Corollary 9.8. $\square$

**Theorem 9.10.** *Suppose that $\zeta_1, \ldots, \zeta_r$ is a filter-regular system of parameters of type $(d_0, \ldots, d_r)$ for $H^*(G, k)$. Then $\zeta_1, \ldots, \zeta_r$ is also filter-regular of type $(d_0, \ldots, d_{r-3}, d_{r-3} - 1, d_{r-3} - 2, d_{r-3} - 2)$.*

*Proof.* Consider the diagram (9.7) with $i = r$, and recall from §2 that $M_{r-1}$ is periodic. So the fact that $\sum_{H \in \mathcal{H}_i} \mathrm{Tr}_{H,G}$ surjects onto $H^j(G, M_{r-1})$ for all large enough values of $j$ means that this holds for all positive values of $j$. Using Theorem 7.1, we see that we may replace $d_{r-2}$ with $d_{r-3} - 1$. Corollary 7.2 now completes the proof of the theorem. $\square$

The following is an improved version of a theorem of Okuyama and Sasaki [22]. Their version has quasi-regular instead of very strongly quasi-regular.

**Corollary 9.11.** *Suppose that the difference between the depth and the Krull dimension of $H^*(G,k)$ is at most two. Then there is a very strongly quasi-regular sequence for $H^*(G,k)$.*

*Proof.* Since the difference between the depth and the Krull dimension (which is equal to $r$) is at most two, we have $a^i_{\mathfrak{m}}(H^*(G,k)) = -\infty$ for $0 \le i \le r-3$ (if there are any such values of $i$). By Theorem 4.5, it follows that there is a filter-regular sequence of type

$$(-1, -2, \ldots, 2-r, d_{r-2}, d_{r-1}, d_r)$$

for some values of $d_{r-2}$, $d_{r-1}$ and $d_r$. Applying Theorem 9.10, it follows that such a sequence is also filter-regular of type

$$(-1, -2, \ldots, 2-r, 1-r, -r, -r),$$

namely a very strongly quasi-regular sequence. □

## 10. Computing group cohomology

In this section, we describe the relationship between the methods of this paper and those developed by Jon Carlson for computing group cohomology [12, 13]. The problem addressed in those papers is as follows. If we try to compute the cohomology of a finite group by beginning a projective resolution of the field of coefficients, then the problem is knowing when to stop. Carlson gives a method for bounding the degrees of the generators and relations of $H^*(G,k)$ in terms of information coming from an initial segment. The method depends on some conjectures in group cohomology, and these conjectures are verified along the way, so that there is no uncertainty in the answer. The conjectures hold in all examples computed so far, including all 2-groups of order at most 64. In this section, we present a method which is guaranteed to work, independently of any conjectures, and which will usually give slightly better bounds than Carlson's. Even a slight improvement in the bound improves the time taken to make the computation by a substantial amount, and makes calculations accessible which were not before.

We begin with the notation. Suppose that we have computed the generators and relations for $H^*(G,k)$ in degrees up to and including $N$. Then we can build an approximation to $H^*(G,k)$ as follows. Choose a minimal set of generators $x_1, \ldots, x_q$ of degrees at most $N$, and let $Q$ be the graded $k$-algebra generated by elements $\tilde{x}_1, \ldots, \tilde{x}_q$, subject only to the relations $\tilde{x}_i \tilde{x}_j = (-1)^{|\tilde{x}_i||\tilde{x}_j|} \tilde{x}_j \tilde{x}_i$. Then there is a homomorphism of graded $k$-algebras $Q \to H^*(G,k)$ taking $\tilde{x}_i$ to $x_i$. Let $J$ be the ideal in $Q$ generated by the elements of the kernel of degree at most $N$. We write $\tau_N H^*(G,k)$ for the quotient $Q/J$, namely the ring defined by the generators and relations of $H^*(G, k)$ of degree at most $N$. Then there is a uniquely defined homomorphism

$$\tau_N H^*(G, k) \to H^*(G, k)$$

taking $\tilde{x}_i + J$ to $x_i$. In general, this homomorphism need not be injective or surjective, but it is an isomorphism up to and including degree $N$. If the degrees of all the generators of $H^*(G, k)$ are at most $N$ then it is surjective, and if the degrees of the relations are also at most $N$ then it is an isomorphism.

**Theorem 10.1.** *Let $r = r_p(G) > 1$. Suppose that $\zeta_1, \ldots, \zeta_r \in \tau_N H^*(G, k)$ form a filter-regular homogeneous set of parameters with $|\zeta_i| = n_i \geq 2$. Suppose further that their images in $H^*(G, k)$ form a homogeneous system of parameters. Set*

$$\alpha = \max_{0 \leq i \leq r-2} \{a_{\mathfrak{m}}^i(\tau_N H^*(G, k)) + i\} \tag{10.2}$$

*($\alpha = -\infty$ if $\tau_N H^*(G, k)$ has depth at least $r - 1$). If*

$$N > \max\{\alpha, 0\} + \sum_{j=1}^{r}(n_j - 1) \tag{10.3}$$

*then the map $\tau_N H^*(G, k) \to H^*(G, k)$ is an isomorphism.*

*Proof.* Let $\mathcal{K}^*(\tau_N H^*(G, k); \zeta_1, \ldots, \zeta_r)$ be the Koszul complex, indexed in the same way as the $_r E_1$ page of the spectral sequence (6.3).

For $0 \leq \ell \leq r$, set

$$d_\ell = \max\{\max_{0 \leq i \leq \ell}\{a_{\mathfrak{m}}^i(\tau_N H^*(G, k)) + i - \ell\}, -\ell - 1\}.$$

Then the sequence $d_0, \ldots, d_r$ satisfies Condition 3.1. Since $a_{\mathfrak{m}}^\ell(\tau_N H^*(G, k)) \leq d_\ell$, it follows from Theorem 4.5 that $\zeta_1, \ldots, \zeta_r$ is a filter-regular sequence of type $(d_0, \ldots, d_r)$ for $\tau_N H^*(G, k)$.

If $0 \leq i \leq r - 2$ then

$$N > \alpha + \sum_{j=1}^{r}(n_j - 1) \geq a_{\mathfrak{m}}^i(\tau_N H^*(G, k)) + i + \sum_{j=1}^{r}(n_j - 1).$$

Since $\max_{0 \leq i \leq \ell}\{d_i + i\} = \max\{\max_{0 \leq i \leq \ell}\{a_{\mathfrak{m}}^i(\tau_N H^*(G, k)) + i\}, -1\}$, setting $s = r - \ell$, we find that for $2 \leq s \leq r$ we have

$$N > d_{r-s} + (r - s) + \sum_{j=1}^{r}(n_j - 1).$$

It follows that $N + s > n_1 + \cdots + n_s + d_{r-s}$, and so by Lemma 3.5, for $2 \leq s \leq r$ and $t \geq N + s$ we have

$$H^{-s,t}(\mathcal{K}^*(\tau_N H^*(G, k); \zeta_1, \ldots, \zeta_r)) = 0. \tag{10.4}$$

There is an obvious map of Koszul complexes

$$\mathcal{K}^*(\tau_N H^*(G, k); \zeta_1, \ldots, \zeta_r) \to {_r E_1^{**}} \tag{10.5}$$

whose target is the $E_1$ page of the spectral sequence ${}_r E_*^{**}$ constructed in §6. Arguing by contradiction, let us suppose that $\tau_N H^*(G,k) \to H^*(G,k)$ is not an isomorphism. Let $N' > N$ be the smallest degree in which $\tau_N H^{N'}(G,k) \to H^{N'}(G,k)$ is not an isomorphism. So there is either a new generator or a new relation in degree $N'$, or both. Since $N' - 2 \geq N - 1 \geq \sum_{j=1}^r (n_j - 1)$, the map (10.5) is an isomorphism on elements of total degree $t - s$ equal to $N' - 2$ or $N' - 1$. Moreover, since $n_i \geq 2$, (10.5) is also an isomorphism for $s > 0$ and $t - s = N'$. And for $t - s$ equal to $N' - 3$, (10.5) is an isomorphism except possibly for $s = r$, $t = -1 + \sum_{j=1}^r n_j$ in the case where $N' = 2 + \sum_{j=1}^r (n_j - 1)$. In this case, the map goes from the zero space to a one dimensional space isomorphic to $\hat{H}^{-1}(G,k)$, namely the element which is sent by $d_r$ to the last survivor, as discussed in §6. In particular, this element is in the kernel of $d_1$.

Taking cohomology with respect to the horizontal differential, we see that

$$H^{-s,t}(\mathcal{K}^*(\tau_N H^*(G,k); \zeta_1, \ldots, \zeta_r)) \to {}_r E_2^{-s,t}$$

is an isomorphism for total degree $t - s$ equal to $N' - 2$ or $N' - 1$, except possibly for $s = 1$, $t = N'$. Combining this with (10.4), we see that for $2 \leq s \leq r$ and $t - s$ equal to $N' - 2$ or $N' - 1$ we have ${}_r E_2^{-s,t} = 0$.

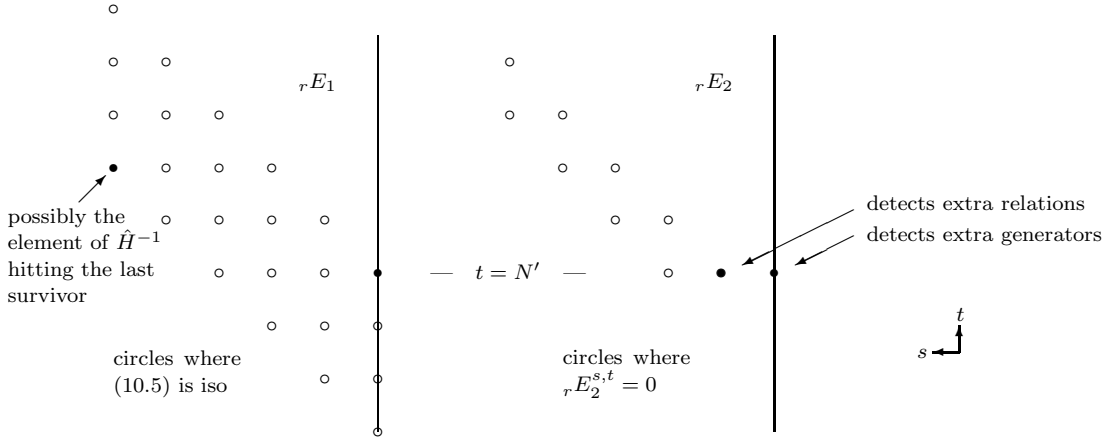

It follows from the case $t - s = N' - 1$ that there is no nonzero element that can hit a new generator in $H^{0,N'}(\mathcal{K}^*(\hat{H}^*(G,k); \zeta_1, \ldots, \zeta_r))$. But the spectral sequence converges to zero, so this is a contradiction. So it follows that $\tau_N H^{N'}(G,k) \to H^{N'}(G,k)$ is surjective. If it is not injective, let $x \in \tau_N H^{N'}(G,k)$ be in the kernel. Then in the Koszul complex, $x$ is the image of some element $y$ of vertical degree $N'$ in $\mathcal{K}^{-1}(\tau_N H^*(G,k); \zeta_1, \ldots, \zeta_r)$. The image of $y$ in ${}_r E_1^{-1,N'}$ survives to give a nonzero element $z \in {}_r E_2^{-1,N'}$, The element $z$ has to be killed somewhere in the spectral sequence. But again, there is no nonzero element in ${}_r E_2^{-s,t}$ with $t - s = N' - 2$ available to hit $z$, giving a contradiction. It follows that $\tau_N$ is an isomorphism as required. □

**Remarks 10.6.** (i) Given elements $\zeta_1, \ldots, \zeta_r \in \tau_N H^*(G, k)$, it is easy to tell whether their images form a homogeneous system of parameters in $H^*(G, k)$. Namely, for a representative of each conjugacy class of maximal elementary abelian $p$-subgroups $E$, we just check that $H^*(E, k)/(\mathrm{res}_{G,E}(\zeta_1), \ldots, \mathrm{res}_{G,E}(\zeta_r))$ has finite length.

(ii) For $N$ large enough, $\tau_N H^*(G, k) \to H^*(G, k)$ is an isomorphism, and the right hand side of the inequality (10.3) is independent of $N$. So for $N$ large enough, (10.3) holds. Since this inequality can be checked just from a knowledge of $\tau_N H^*(G, k)$, the theorem gives a computable bound on the degrees of the generators and relations of $H^*(G, k)$, without referring to any conjectures.

(iii) Computing $\alpha$ in practise can be accomplished using Theorem 4.5.

(iv) The definition (10.2) of $\alpha$ implies that $\alpha \leq \mathsf{Reg}\,\tau_N H^*(G, k)$. So if Conjecture 1.1 holds, then we only ever have to compute up to degree
$$N = 1 + \sum_{j=1}^{r}(n_j - 1).$$
In particular, it follows from Corollary 9.11 that this is true if the depth and Krull dimension of $H^*(G, k)$ differ by at most two.

(v) If the Sylow $p$-subgroup of $G$ has rank at least two, then by Duflot's theorem [16], the depth of $H^*(G, k)$ is at least two. Then a theorem from [7] shows that there are no nonzero products from negative degrees to positive degrees in Tate cohomology. It follows that Theorem 10.1 can be improved in this case by replacing (10.3) by the inequality
$$N \geq \max\{\alpha, 0\} + \sum_{j=1}^{r}(n_j - 1).$$
In the proof of the theorem, the diagram for $_rE_1$ becomes

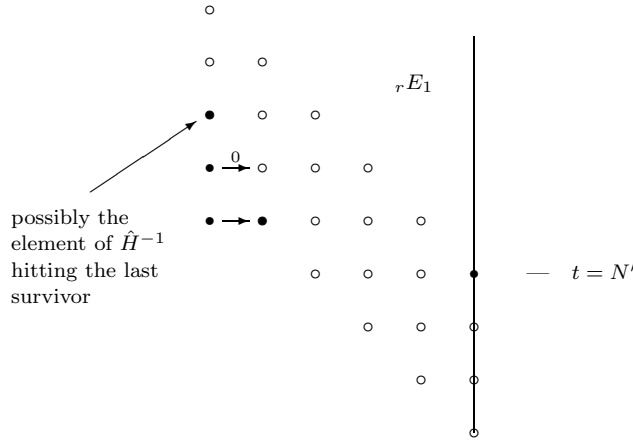

It is possible that there are unwanted elements of $\hat{H}^{-1}(G,k)$ appearing in the second column from the left, but they are hit from the leftmost column, so they do not damage the argument.

(vi) Dave Rusin [23] computed the cohomology of all the 2-groups of order up to 32, using the Eilenberg–Moore spectral sequence. This method becomes impractical for the groups of order 64. Jon Carlson has computed the cohomology of the groups of order 64 by computing the beginning of a resolution explicitly using the computer algebra package MAGMA, and using the methods of [13] to know when to stop computing. It turns out that for all the 2-groups of order at most 64, the depth and the Krull dimension differ by at most two, so that $\mathsf{Reg}\, H^*(G,k)$ is guaranteed to be zero. In just one case, namely the Sylow 2-subgroup of $PSL(3,4)$, $H_{\mathfrak{m}}^{r,-r} H^*(G,k)$ turned out to have dimension bigger than one. So the last survivor does not span this two dimensional space in this case.

## References

1. D. J. Benson, *Representations and Cohomology II: Cohomology of groups and modules*, Cambridge Studies in Advanced Mathematics, vol. 31, Cambridge University Press, 1991, reprinted in paperback, 1998.
2. \_\_\_\_\_\_, *Polynomial invariants of finite groups*, London Math. Soc. Lecture Note Series, vol. 190, Cambridge University Press, 1993.
3. \_\_\_\_\_\_, *The image of the transfer*, Archiv der Math. (Basel) **61** (1993), 7–11.
4. \_\_\_\_\_\_, *Conway's group $Co_3$ and the Dickson invariants*, Manuscripta Math. **85** (1994), 177–193.
5. \_\_\_\_\_\_, *Modules with injective cohomology, and local duality for a finite group*, New York Journal of Mathematics **7** (2001), 201–215.
6. \_\_\_\_\_\_, *Commutative algebra in the cohomology of groups*, Lecture notes for three lectures given at the MSRI Commutative Algebra opening conference in September, 2002.
7. D. J. Benson and J. F. Carlson, *Products in negative cohomology*, J. Pure & Applied Algebra **82** (1992), 107–129.
8. \_\_\_\_\_\_, *Projective resolutions and Poincaré duality complexes*, Trans. Amer. Math. Soc. **132** (1994), 447–488.
9. D. J. Benson and C. W. Wilkerson, *Finite simple groups and Dickson invariants*, Homotopy theory and its applications, Cocoyoc, Mexico, 1993, Contemp. Math., vol. 188, 1995, pp. 39–50.
10. D. Bourguiba and S. Zarati, *Depth and the Steenrod algebra, with an appendix by J. Lannes*, Invent. Math. **128** (1999), 589–602.
11. M. P. Brodmann and R. Y. Sharp, *Local cohomology*, Cambridge Studies in Advanced Mathematics, vol. 60, Cambridge University Press, 1998.
12. J. F. Carlson, *Problems in the calculation of group cohomology*, Proceedings of the Euroconference on Computational Methods for Representations of Groups and Algebras, Progress in Mathematics, vol. 173, 1999, pp. 107–120.

13. ______, *Calculating group cohomology: Tests for completion*, J. Symb. Comp. **31** (2001), 229–242.
14. J. F. Carlson, C. Peng, and W. W. Wheeler, *Transfer maps and virtual projectivity*, J. Algebra **204** (1998), 286–311.
15. H. Cartan and S. Eilenberg, *Homological algebra*, Princeton Mathematical Series, no. 19, Princeton Univ. Press, 1956.
16. J. Duflot, *Depth and equivariant cohomology*, Comment. Math. Helvetici **56** (1981), 617–637.
17. D. Eisenbud, *Commutative algebra, with a view towards algebraic geometry*, Graduate Texts in Mathematics, vol. 150, Springer-Verlag, Berlin/New York, 1995.
18. D. Eisenbud and S. Goto, *Linear free resolutions and minimal multiplicity*, J. Algebra **88** (1984), 89–133.
19. L. Evens, *The cohomology ring of a finite group*, Trans. Amer. Math. Soc. **101** (1961), 224–239.
20. ______, *Cohomology of groups*, Oxford University Press, 1991.
21. J. P. C. Greenlees, *Commutative algebra in group cohomology*, J. Pure & Applied Algebra **98** (1995), 151–162.
22. T. Okuyama and H. Sasaki, *Private communication*, Chiba, Japan, 2000.
23. D. J. Rusin, *The cohomology of groups of order* 32, Math. Comp. **53** (1989), 359–385.
24. N. V. Trung, *Reduction exponent and degree bound for the defining equations of graded rings*, Proc. Amer. Math. Soc. **101** (1987), 229–236.
25. ______, *The largest non-vanishing degree of graded local cohomology modules*, J. Algebra **215** (1999), 481–499.
26. C. W. Wilkerson, *A primer on Dickson invariants*, Proc. of the Northwestern Homotopy Theory Conf., Contemp. Math., vol. 19, American Math. Society, 1983, pp. 421–434.

DEPARTMENT OF MATHEMATICS, UNIVERSITY OF GEORGIA, ATHENS GA 30602, USA
*E-mail address*: djb@byrd.math.uga.edu